\documentclass[10pt,a4paper,twoside]{amsart}
\RequirePackage[l2tabu, orthodox]{nag}

\usepackage[utf8]{inputenc}
\usepackage[english]{babel}
\usepackage[T1]{fontenc}
\usepackage{microtype, fancyhdr, lmodern, xcolor}

\usepackage[textsize=9]{todonotes}

\usepackage{amsfonts, amsmath, amsthm, amssymb, mathrsfs, amscd}		
\numberwithin{equation}{section}	
\usepackage{faktor}
\allowdisplaybreaks

\usepackage{url, enumitem}
\usepackage[noadjust]{cite}

\usepackage{multirow,bigdelim, caption}
\usepackage{booktabs}

\usepackage{graphicx}

\usepackage{hyperref}

\theoremstyle{plain}
\newtheorem{thm}{Theorem}[section]
\newtheorem{prop}[thm]{Proposition}

\newtheorem{lem}[thm]{Lemma}

\theoremstyle{remark}
\newtheorem{rmk}[thm]{Remark}

\newtheorem{step}{Step}

\theoremstyle{definition}

\newcommand\Zz{\mathbb{Z}}

\newcommand\Rr{\mathbb{R}}

\DeclareMathOperator{\Aut}{Aut}	

\newcommand\LBE[1]{LB_{#1}^{ext}}	

\newcommand \F[1]{F_{#1}}

\DeclareMathOperator{\M}{\mathcal{M}}	

\newcommand\ie{{\textit{i.e.}}}
\newcommand\ii{i}
\newcommand\inv{^{-1}}

\newcommand\jj{j}

\newcommand\nn{n}
\newcommand\nno{{n-1}}

\makeatletter
\newcommand{\oset}[2]{%
  {\mathop{#2}\limits^{\vbox to -.5\ex@{\kern-\tw@\ex@
   \hbox{\scriptsize #1}\vss}}}}
\makeatother

\newcommand\R[1]{R_{#1}}	

\newcommand\rr[1]{\rho_{\hspace{-0.3ex}#1}^{\null}}	
\newcommand\rrq[2]{\rho_{\hspace{-0.3ex}#1}^{#2}}	

\newcommand\sig[1]{\sigma_{\hspace{-0.3ex}#1}^{\null}} 
\newcommand\sigg[2]{\sigma_{\hspace{-0.3ex}#1}^{#2}}	
\newcommand\siginv[1]{\sigg{#1}{-1}}	

\DeclareMathOperator{\Ring}{\mathcal{R}}	
\DeclareMathOperator{\Link}{\mathcal{L}}	

\begin{document}

\title{On the group of ring motions of an H-trivial link}

\author[Damiani]{Celeste Damiani}
\address{Department of Mathematics,
Osaka City University,
Sugimoto, Sumiyoshi-ku,
Osaka 558-8585, Japan}
\email{celeste.damiani@math.cnrs.fr}

\author[Kamada]{Seiichi Kamada}
\address{Department of Mathematics,
Osaka City University,
Sugimoto, Sumiyoshi-ku,
Osaka 558-8585, Japan}
\email{skamada@sci.osaka-cu.ac.jp}



\date{\today}

\begin{abstract}

In this paper we compute a presentation for the group of ring motions 
of the split union of a Hopf link with Euclidean components and a Euclidean circle.
A key part of this work is the study of a short exact sequence 
of groups of ring motions of general ring links in~$\Rr^3$. This sequence allowed us to 
build the main result from the previously known case of the ring group with one component, which 
a particular case of the ring groups studied by Brendle and Hatcher. 
This work is a first step towards the computation of a presentation for 
groups of motions of H-trivial links with an arbitrary number of components. 
\end{abstract}

\maketitle


\section{Introduction}
\label{S:Intro}

An \emph{H-trivial link} of type $(m,n)$ is a link in $\Rr^3$ 
which is ambiently isotopic to the split union of $m$ Hopf links and $n$ trivial knots.  
When~$m=0$, it is a trivial link with $n$ components. 
H-trivial links are a generalization of trivial links, 
and play an important role in normal forms of 
immersed surface-links in~$\Rr^4$~\cite{Kamada-Kawamura:2017,  Kamada-Kawauchi-Kim-Lee:2017}.

A \emph{ring} in $\Rr^3$ is a circle in the strict Euclidian sense, \ie, a round circle on a plane in~$\Rr^3$.  
We call a link in $\Rr^3$ a \emph{ring link} if each component is a ring. 
The \emph{ring group} $R_n$ (of a trivial ring link with $n$ components) 
was introduced by Brendle and Hatcher~\cite{BrendleHatcher:2013} 
as the fundamental group of the space 
of all configurations of ring links 
which are equivalent, as ring links\footnote{
The original definition of $R_n$ in~\cite{BrendleHatcher:2013} is the 
fundamental group of the space 
of all configurations of ring links 
which are equivalent \emph{as links} to a trivial link with $n$ components. 
If a ring link is equivalent as a link to a  trivial ring link, then 
it is equivalent as a ring link. This fact is asserted in~\cite{BrendleHatcher:2013}.   
}, to a trivial ring link with $n$ components.  
We generalize this notion to the ring group $R_{m,n}$ 
as the fundamental group of the space of all configurations of 
ring links 
which are equivalent, as ring links, to an H-trivial ring link of type~$(m,n)$.  
We give presentations of the ring groups $R_{m,n}$ for~$(m,n)= (0,1)$,~$(1,0)$, and~$(1,1)$.  
Some basic properties of the group $R_{m,n}$ are also given.

The paper is structured as follows: in Section~\ref{S:generalities_motions} we give the basic definitions 
concerning ring motions, and we discuss some tools and properties. 
In Section~\ref{S:Motion} we review known results about the ring group $\R\nn$ of a trivial link, discussing its relation
with the motion group of a trivial link studied in~\cite{Dahm:Thesis} and~\cite{Goldsmith:MotionGroupsTrivial},
and recalling a 
complete presentation given in~\cite{BrendleHatcher:2013} (Proposition~\ref{P:BrendleHatcher}). In Section~\ref{S:Ring}
we introduce an exact sequence for groups of ring motions of ring links (Proposition~\ref{P:ExactSequenceRing}) on which we rely to 
find presentations for many of the considered groups.
In Section~\ref{S:Circle} we focus on the particular case
of the ring group $\R1$ of just one ring. Here we give an alternative argument for the proof of its presentation (Lemma~\ref{L:RingCircle}). This serves 
as a strategy model for the case of the ring group $R_{1,0}$ of a Hopf link, treated in Section~\ref{S:Hopf} (Lemma~\ref{L:PresentationHopf}).
Finally in Section~\ref{S:HopfCircle} we join all preliminary results, and using standard techniques to 
write presentations for group extensions we give a presentation for the group 
of ring motions  $R_{1,1}$ of a H-trivial ring link of type~$(1, 1)$ in the main result of this paper~(Theorem~\ref{T:HopfCirclepresentation}).

\section{Ring motions and motions of links} 
\label{S:generalities_motions}

Let $M$ be a $3$-manifold in~$\Rr^3$.  
A link in $M$ is called a \emph{ring link} if each component is a ring. 
Two ring links $L$ and $L'$ in $M$ are \emph{equivalent} (as ring links in $M$) if there exists an isotopy 
$\{ L_t \}_{t \in [0,1]}$ through ring links $L_t$ in $M$ with $L_0=L$ and~$L_1=L'$.  

For a ring link $L$ in~$M$, let 
$\Ring(M, L)$ be the space of all 
configurations of 
ring links 
which are equivalent, as ring links in $M$, to $L$.  This space has $L$ as base point. 
The \emph{ring group} of $L$ in~$M$, 
denoted by~$R(M,L)$, is 
the fundamental group~$\pi_1(\Ring(M, L))$.  

Let $L_{m,n}$ be a ring link in $\Rr^3$ which is a \emph{split} union 
of $m$ Hopf links and $n$ trivial knots, namely, each Hopf link (and each trivial knot component) can be separated from the other by a convex hull in~$\Rr^3$.  
The \emph{ring group} $R_{m,n}$ is 
the ring group $R(\Rr^3, L_{m,n})$ of~$L_{m,n}$, \ie, the fundamental group of the space of all configurations of ring links which are equivalent, as ring links, to ~$L_{m,n}$.
This group does not depend on the choice of a base point~$L_{m,n}$.


A \emph{ring motion} of a ring link $L$ in $M$ is a loop in the based space~$\Ring(M, L)$, 
which is presented by a  
$1$-parameter family $\{L_t\}_{t \in [0,1]}$ of ring links in $M$ 
with~$L =L_0 = L_1$.  
The \emph{stationary motion} or the \emph{trivial motion} of $L$ is a ring motion $\{L_t\}_{t \in [0,1]}$ with $L=L_t$ for all $t \in [0,1]$. 
Two ring motions are said to be 
\emph{equivalent} (as ring motions) or \emph{homotopic}
if they are homotopic through ring motions of $L$ in~$M$.  
The product of two ring motions are defined by concatenation. 
The set of equivalence classes of ring motions of $L$ in $M$ forms a group. This is, by definition,   
the ring group~$R(M, L)$.

Ring groups are related to motion groups
as introduced by
Dahm~\cite{Dahm:Thesis} and Goldsmith~\cite{Goldsmith:MotionGroupsTrivial, Goldsmith:MotionGroupsTorus}.  
Let $M$ be a $3$-manifold and $L$ a link in~$M$. 
Roughly speaking, a \emph{motion} of $L$ in $M$ 
is a $1$-parameter family $\{L_t\}_{t \in [0,1]}$ of links in $M$ 
with $L =L_0 = L_1$ such that there exists 
an ambient isotopy $\{f_t\}_{t \in [0,1]}$ of $M$ with compact support 
and 
such that $L_t = f_t(L)$ for~$t \in [0,1]$.  
Two motions are said to be \emph{equivalent} (as motions) 
or \emph{homotopic} 
if they are homotopic through motions of $L$ in~$M$.  
The product of two motions is defined by concatenation.  
The set of equivalence classes of motions of $L$ in $M$ 
forms a group, which is the \emph{motion group} of $L$ in $M$ 
and is denoted by~$\M(M, L)$.  
For 
a detailed treatment of motions and motion groups,
we refer to Dahm~\cite{Dahm:Thesis} and Goldsmith~\cite{Goldsmith:MotionGroupsTrivial, Goldsmith:MotionGroupsTorus}.

For a ring link $L$ in a $3$-manifold~$M \subset \Rr^3$, 
there is a natural homomorphism 
\[
R(M, L) \to \M(M, L). 
\]
This map is an isomorphism when $M= \Rr^3$ and $L$ is a trivial ring link~{\cite[Theorem 1]{BrendleHatcher:2013}}, 
{\cite[Theorem~3.10]{Damiani:Journey}}.

\section{The ring group and the motion group of a trivial link}
\label{S:Motion}

In this section we recall 
some known results about
the group $R_{0, n}=  R_n = R(\Rr^3, L) \cong \M(\Rr^3, L)$ 
of a trivial link $L$ with $n$ components.  

Let $L$ be a link in~$\Rr^3$. 
The \emph{Dahm homomorphism} is a well-defined homomorphism  
\[ 
D \colon \M( \Rr^3 , L)  \longrightarrow \Aut(\pi_1(\Rr^3 \setminus L)),   
\]
defined as follows.
Let $\{ L_t \}_{t \in [0,1]}$ be a motion of $L$ in~$\Rr^3$, and $p$ a base point far from the motion. 
Let $A  \subset \Rr^3 \times [0,1]$ be the annulus with $A \cap \Rr^3 \times \{t\} = L_t \times \{t\}$ for~$t \in [0,1]$.  
Consider the 
automorphism
$(i_1)^{-1}_\ast \circ (i_0)_\ast \colon \pi_1(\Rr^3 \setminus L; p) \to \pi_1(\Rr^3 \setminus L; p)$, 
where~$i_k$, for $(k=0,1)$, is the inclusion map of $\Rr^3 \setminus L = (\Rr^3 \setminus L) \times \{k\}$ to~$\Rr^3 \times [0,1] \setminus A$.  
Then $D(\{ L_t \}_{t \in [0,1]})$  is defined by this automorphism. 

The Dahm homomorphism is also defined on the ring group $R(\Rr^3, L)$,   
\[ D \colon R( \Rr^3 , L)  \longrightarrow \Aut(\pi_1(\Rr^3 \setminus L)).    
\]

Let~$\nn \geq 1$, and  $C = C_1 \sqcup \cdots \sqcup C_\nn$ 
be a trivial (ring) link with $\nn$ components in~$\Rr^3$, 
with $C_i =\{ (x,y,0) \in \Rr^3 \mid (x - i)^2 + y^2 = (1/4)^2 \}$ for~$i=1, \dots, \nn$.  

The fundamental group $\pi_1(\Rr^3 \setminus C)$ is the free group $\F\nn$ of rank~$\nn$ generated by $x_1, \dots, x_{\nn}$, where $x_i$ is the element represented by a positively oriented meridian loop of $C_k$ with respect to the counterclockwise orientation of~$C_k$.    

The two following results display some basic properties 
for the motion group $\M(\Rr^3, C)$ and the ring group $R_n$. 
These will lead to explaining the relation between the two, and to recalling a presentation for these groups.

\begin{thm}[{\cite[Theorems~5.3 and 5.4]{Goldsmith:MotionGroupsTrivial}}]
\label{T:Gold} 
\mbox{}
\begin{enumerate}[label={(\arabic*)}]
\item The Dahm homomorphism  
\[ 
D \colon \M(\Rr^3, C) \longrightarrow \Aut(\F\nn)
 \] 
is injective.  

\item
The motion group $\M(\Rr^3, C)$ is generated by the following types of motions: 
\begin{itemize}
\item Permute the $i$th and the $(i+1)$st rings by pulling the $i$th ring through the $(i+1)$st ring.  
\item Permute the $i$th and the $(i+1)$st rings by passing the $i$th ring around the $(i+1)$st ring.  
\item Reverse the orientation of the $i$th ring by rotationg it by 180 degrees around the $x$-axis. 
\end{itemize}

\item The above generators correspond to the following automorphisms of~$F_\nn$:  
\begin{align}
\label{E:sigma}
\sig\ii &: \begin{cases}
            x_\ii \mapsto x_{\ii+1}; &\\
            x_{\ii+1} \mapsto x_{\ii+1}\inv x_\ii x_{\ii+1}; &\\
            x_\jj \mapsto x_\jj, \ &\text{for} \ \jj \neq \ii, \ii+1.
        \end{cases} \\
\label{E:rho}
\rr\ii &: \begin{cases}
            x_\ii \mapsto x_{\ii+1}; & \\
            x_{\ii+1} \mapsto x_\ii; &\\
            x_\jj \mapsto x_\jj, \ &\text{for} \ \jj \neq \ii, \ii+1.
		\end{cases} \\
\label{E:tau}
\tau_\ii &: \begin{cases}
            x_\ii \mapsto x\inv_\ii; &\\
            x_\jj \mapsto x_\jj, \ &\text{for} \ \jj \neq \ii.
        \end{cases}
\end{align}

\item The image of the Dahm homomorphism, \ie, the subgroup of $\Aut(F_\nn)$ generated by the above automorphisms,
is the
group of automorphisms of $F_\nn$ of the form $\alpha \colon x_\ii \mapsto w_\ii\inv x^{\pm 1}_{\pi({\ii})} w_\ii$, 
where $\pi$ is a permutation of the indices and $w_\ii$ is a word in~$F_\nn$ (compare with the \emph{group of conjugating automorphisms}~\cite{Savushkina:1996}, also known as \emph{group of permutation-conjugacy automorphisms}~\cite{Suciu-Wang:2017}). 
\end{enumerate}

\end{thm}

\begin{thm}[{\cite[Theorem~1]{BrendleHatcher:2013}}]
\label{T:RelaxingCircles}
Let $\Ring_n$ be the 
configuration space of
ring links which are equivalent to $C$ and 
let $\Link_n$  be the space of all smooth links equivalent to $C$.  The inclusion of $\Ring_n$ 
into  $\Link_n$ is a homotopy equivalence.
\end{thm} 

Leaning on Theorem~\ref{T:RelaxingCircles} it is possible to show
that there is a natural isomorphism between
$R_n = R( \Rr^3, C) $ and $\M( \Rr^3 , C)$~\cite[Theorem~3.10]{Damiani:Journey}.
Thus the 
statement of Theorem~\ref{T:Gold} holds for the ring group~$R_n$ too. 

\begin{rmk} Our notations $\sig\ii, \rr\ii, \tau_\ii$ for the motions and the automorphisms 
in Theorem~\ref{T:Gold} are different from those used in \cite{Goldsmith:MotionGroupsTrivial} or~
\cite{BrendleHatcher:2013}. However they coincide with the ones used in~\cite{Damiani:Journey}, 
where this group is called \emph{extended loop braid group~$\LBE\nn$}.
\end{rmk}

\begin{prop}[{\cite[Theorem~3.7]{BrendleHatcher:2013}}]
\label{P:BrendleHatcher}
The group $\R\nn$ admits a presentation given by the sets of generators $\{\sig\ii, \rr\ii \mid \ii=1, \dots , \nno \}$ and $\{\tau_\ii \mid \ii=1, \dots , \nn \}$ subject to relations:
\begin{equation}
\label{E:presentation}
\begin{cases}
\sig{i} \sig j = \sig j \sig{i}  \, &\text{for } \vert  i-j\vert > 1\\
\sig{i} \sig {i+1} \sig{i} = \sig{i+1} \sig{i} \sig{i+1} \, &\text{for } i=1, \dots, \nn-2 \\
\rr{i} \rr j = \rr j \rr{i}  \, &\text{for }  \vert  i-j\vert > 1\\
\rr\ii\rr{i+1}\rr\ii = \rr{i+1}\rr\ii\rr{i+1} \, &\text{for }  i=1, \dots, \nn-2  \\
\rrq{i}2 =1 \, &\text{for }  i=1, \dots, \nno \\
\rr{i} \sig{j} = \sig{j} \rr{i}   \, &\text{for }  \vert  i-j\vert > 1\\
\rr{i+1} \rr{i} \sig{i+1} = \sig{i} \rr{i+1} \rr{i} \,  &\text{for }  i=1, \dots, \nn-2  \\
\sig{i+1} \sig{i} \rr{i+1} = \rr{i} \sig{i+1} \sig{i}  \,  &\text{for }  i=1, \dots, \nn-2 \\
\tau_{i} \tau_j = \tau_j \tau_{i}  \, &\text{for }    i \neq j \\
\tau_\ii^2=1 \, &\text{for }  i=1, \dots, \nn \\
\sig{i} \tau_j = \tau_j \sig{i}  \, &\text{for }  \vert  i-j\vert > 1 \\
\rr{i} \tau_j = \tau_j \rr{i}  \, &\text{for }  \vert  i-j\vert > 1 \\
\tau_i \rr\ii = \rr\ii \tau_{i+1} \, &\text{for }  i=1, \dots, \nno  \\
\tau_\ii \sig\ii = \sig\ii \tau_{i+1}  \, &\text{for } i=1, \dots, \nno \\
\tau_{i+1} \sig\ii = \rr\ii \siginv\ii \rr\ii \tau_\ii  \, &\text{for } i=1, \dots, \nno. \\
\end{cases}
\end{equation}
\end{prop}

\section{Extensions and projections}
\label{S:Ring}

Let $L_1$ and $L_2$  be ring links 
in a $3$-manifold $M \subset \Rr^3$ 
with~$L_1 \cap L_2 = \emptyset$.  

We say that 
 a ring motion $ \{ L_{1(t)} \}_{t \in [0,1]}$ of $L_1$ in $M$ and 
a ring motion  $ \{ L_{2(t)} \}_{t \in [0,1]}$ of $L_2$ in $M$ 
are \emph{disjoint}  if $L_{1(t)} \cap L_{2(t)} = \emptyset$ for all~$t \in [0,1]$. 
In this case,  $\{ L_{1(t)}  \sqcup  L_{2(t)} \}_{t \in [0,1]}$ is a ring motion of $L_1 \sqcup L_2$ in~$M$.  
We denote this ring motion by $\{ L_{1(t)} \}_{t \in [0,1]} \sqcup \{ L_{2(t)} \}_{t \in [0,1]}$ 
and call it the \emph{union} of the motions    $ \{ L_{1(t)} \}_{t \in [0,1]}$ and~$ \{ L_{2(t)} \}_{t \in [0,1]}$.  
 
We denote by $R(M, L_1, L_2)$ the subgroup of the ring group $R(M, L_1 \sqcup L_2)$ 
consisting of equivalence classes of ring motions 
which can be written as the union of a motion of $L_1$ and a motion of~$L_2$.  
It is a subgroup of index two if and only if there exists a ring motion of $L_1 \sqcup L_2$ in $M$ 
which interchanges $L_1$ and~$L_2$.  
Otherwise, $R(M, L_1, L_2) = R(M, L_1 \sqcup L_2)$.  

For a ring motion $\{ L_{2(t)}\}_{t \in [0,1]}$ of $L_2$ in~$M \setminus L_1$, 
we have a ring motion 
$\{ L_1 \}_{t \in [0,1]} \sqcup \{ L_{2(t)}\}_{t \in [0,1]}$ 
of $L_1 \sqcup L_2$ in~$M$, 
where   $\{ L_1 \}_{t \in [0,1]}$ is the stationary motion of~$L_1$. 
We call it the \emph{extension} of $\{ L_{2(t)}\}_{t \in [0,1]}$ with~$L_1$, 
and we denote it by~$e(\{ L_{2(t)}\}_{t \in [0,1]})$.  
We have a well-defined homomorphism 
\[ e \colon R(M \setminus L_1, L_2) \longrightarrow R(M, L_1, L_2) \]
with $e([ \{ L_{2(t)}\}_{t \in [0,1]} ] ) = [e(\{ L_{2(t)}\}_{t \in [0,1]}) ]$.   

Let 
\[
p_1 \colon R(M, L_1, L_2) \longrightarrow R(M, L_1)
\]
be the homomorphism sending 
$[ \{ L_{1(t)} \}_{t \in [0,1]} \sqcup \{ L_{2(t)} \}_{t \in [0,1]} ]$ to 
$ [ \{ L_{1(t)} \}_{t \in [0,1]} ]$.  

\begin{prop}
\label{P:ExactSequenceRing}
Let $L_1$ and $L_2$ be disjoint ring links in~$M \subset \Rr^3$. 
Consider the composition of $e$ and~$p_1$:
\begin{equation}
\label{E:epRing}
 \begin{CD}
R(M \setminus L_1, L_2)  @>e>>  R(M, L_1, L_2) @>p_1>> R(M, L_1).
\end{CD} 
\end{equation} 
Then 
${\mathrm Im}~{e} \subset {\mathrm Ker}~{p_1}$.  
\end{prop}

\proof
Follows directly from the definitions
of the applications $e$ and~$p_1$.  
\endproof

Although it appears that sequence~\eqref{E:epRing} is exact in many cases, 
few examples are known to the authors at this moment.  For example, in the case 
of trivial links due to 
\cite{BrendleHatcher:2013} and in the cases which we discuss in this paper, 
sequence~\eqref{E:epRing} is exact. 

\begin{rmk}[]
\label{P:ExactSequence} 
Let $L_1$ and $L_2$ be disjoint links in a $3$-manifold~$M$. It is known~{\cite[Proposition~3.10]{Goldsmith:MotionGroupsTrivial}}
that the following sequence on motion groups is exact:
\[ \begin{CD}
\M(M \setminus L_1, L_2)  @>e>> \M(M, L_1, L_2) @>p_1>> \M(M, L_1).
\end{CD} \]
\end{rmk}

\section{The ring group of a ring}
\label{S:Circle}

First we observe the ring group of a ring $C$ in~$\Rr^3$. 
Let $C$ be the unit ring~$\{ (x, y, 0) \in \Rr^3 \mid x^2 + y^2 =1 \}$. 
In \cite{Goldsmith:MotionGroupsTrivial} and  
\cite{BrendleHatcher:2013} it is shown that the ring group $R(\Rr^3, C)$ and 
the motion group $\M(\Rr^3, C)$ 
are cyclic groups of order $2$ 
generated by the class of a ring motion of $C$ rotating it $180$ degrees about   the $y$-axis.   

Let $R_x(\varphi), R_y(\varphi), R_z(\varphi)$  denote (counterclockwise) rotations of $\Rr^3$ about the $x$-axis, the $y$-axis and the $z$-axis by angle $\varphi$. 
These are identified with special orthogonal matrices as follows: 
\[
R_x(\varphi)= 
\left( 
\begin{array}{ccc} 
1 & 0 & 0 \\ 
0 & \cos \varphi & -\sin \varphi \\ 
0 & \sin \varphi & \cos \varphi 
\end{array}
\right), 
\quad 
R_y(\varphi)= 
\left( 
\begin{array}{ccc} 
\cos \varphi & 0 & \sin \varphi \\ 
0 & 1 & 0 \\ 
-\sin \varphi & 0 & \cos \varphi 
\end{array}
\right) 
\]

\[
R_z(\varphi)= 
\left( 
\begin{array}{ccc} 
\cos \varphi & -\sin \varphi & 0 \\ 
\sin \varphi & \cos \varphi & 0 \\ 
0 & 0 & 1 
\end{array}
\right). 
\] 

Let $\tau_C $ be the element of $R(\Rr^3, C)$ represented by a ring motion $\{ R_y(\pi t)(C) \}_{t \in [0,1]}$, \ie, the $180$ degrees rotation about the $y$-axis.  

\begin{lem}[{\cite[Theorem~3.7]{BrendleHatcher:2013}}, {\cite[Theorem~5.3]{Goldsmith:MotionGroupsTrivial}}]
\label{L:RingCircle}
The ring group $R(\Rr^3, C)$, which is isomorphic  to the motion group $\M(\Rr^3, C)$, admits the presentation 
\[
\langle \tau_C \mid \tau_C^2 = 1 \rangle. 
\]
\end{lem}

\proof 
We only show the case of~$R(\Rr^3, C)$.  
Any ring $L$ in $\Rr^3$ is determined uniquely by the position of 
the center~$c(L) \in \Rr^3$, 
the radius~$r(C) \in \Rr_{>0}$,
and an element $g(C)$ of the Grassman manifold $G(2,3)$ 
of unoriented $2$-planes through the origin $O$ in~$\Rr^3$,
which is obtained from the plane $H(C)$ containing $C$  
by sliding it along the vector~$\overrightarrow{c(L)O}$.  
Thus the space of rings in $\Rr^3$ is identifies with 
$\Rr^3 \times \Rr_{>0} \times G(2,3)$. 
There is a deformation retract to the subspace~$\{O\} \times \{1\} \times G(2,3) \cong G(2,3)$.  
The fundamental group of $G(2,3)$ is a cyclic group of order $2$ 
generated by a loop which rotates the $xy$-plane by $180$ degrees about the $y$-axis.  
This corresponds~$\tau_C \in R(\Rr^3, C)$.  
\endproof 

The proof above suggests a strategy to deform a ring motion
to a \lq\lq standard\rq\rq~ ring motion, 
which is used later for the case of a Hopf link.

\section{The ring group of a Hopf link}
\label{S:Hopf}

Let $H_1$ and $H_2$ be unit rings in $\Rr^3$ with 
$H_1 = \{ (x, y, 0) \in \Rr^3 \mid x^2 + y^2 =1 \}$ and~$H_2 = \{ (0, y, z) \in \Rr^3 \mid (y-1)^2 + z^2 =1 \}$.   
The \emph{positive standard rotation} of $H_2$ along $H_1$ is 
a ring motion $\{ R_z(2\pi t)(H_2) \}_{t \in [0,1]}$ of $H_2$ in $\Rr^3 \setminus H_1$ or in~$\Rr^3$, 
and the \emph{negative standard rotation} of $H_2$ along $H_1$ is a ring motion~$\{ R_z(-2\pi t)(H_2) \}_{t \in [0,1]}$. 

\begin{lem}
\label{L:MH2}
The ring group $R(\Rr^3 \setminus H_1, H_2)$ admits the presentation 
\[
\langle \ell \mid ~~~  \rangle,  
\]
where $\ell$ is represented by the positive standard rotation of $H_2$ along~$H_1$. 
\end{lem}

First we introduce the \emph{rotation number} of a ring motion of $H_2$ in $\Rr^3 \setminus H_1$ such that we obtain a 
homomorphism ${\mathrm rot}: R(\Rr^3 \setminus H_1, H_2) \to \Zz$ with~${\mathrm rot}(\ell)=1$. 


Given an orientation on~$H_2$, note that $H_2$ always comes back to itself with the same orientation after any ring motion $H_2$ in~$\Rr^3 \setminus H_1$.
Thus, a ring motion of $H_2$ in $\Rr^3 \setminus H_1$   
induces a continuous map from~$H_2 \times S^1 \to \Rr^3 \setminus H_1$,
and hence a homomorphism $H_2(H_2 \times S^1; \Zz) \to H_2(\Rr^3 \setminus H_1; \Zz)$
on the $2$nd homology groups.  
We call it the homomorphism on the $2$nd homology groups induced from the motion of~$H_2$. 
Note that if two ring motions are homotopic as ring motions then their induced homomorphisms are the same.    

Note that $H_2(H_2 \times S^1; \Zz) \cong \Zz$ and $H_2(\Rr^3 \setminus H_1; \Zz) \cong \Zz$ and the 
homomorphism induced from the positive standard rotation of $H_2$ along $H_1$ sends a generator to a generator. 
Choose generators of $H_2(H_2 \times S^1; \Zz) \cong \Zz$
and $H_2(\Rr^3 \setminus H_1; \Zz) \cong \Zz$ so that 
the homomorphism induced from the positive standard rotation of $H_2$ along $H_1$ 
sends $1 \in \Zz$ to~$1 \in \Zz$. 
 The \emph{rotation number} of the motion is the integer 
 which is the image of $1$ under the induced homomorphism 
 on the $2$nd homology groups. 
 It yields the desired homomorphism ${\mathrm rot}: R(\Rr^3 \setminus H_1, H_2) \to \Zz$ with~${\mathrm rot}(\ell)=1$. 

We call a ring motion $\{ L_t \}_{t \in [0,1]}$ of $H_2$ 
in $\Rr^3 \setminus H_1$ a \emph{normal} ring motion if 
there is a continuous map $\phi: [0,1] \to \Rr$ 
such that $L_t= R_z(2 \pi \phi(t))(H_2)$ for all~$t \in [0,1]$.  
For a normal ring motion, 
$\phi(1) - \phi(0) \in \Zz$ is the rotation number.  
We have that 
the equivalence class of a normal ring 
motion is~$\ell^{\phi(1) - \phi(0)} \in R(\Rr^3 \setminus H_1, H_2)$.    

\proof[Proof of Lemma~\ref{L:MH2}] 
It is sufficient to show that $R(\Rr^3 \setminus H_1, H_2)$ is generated by~$\ell$,
by using the rotation number.  
Let $\{ L_t \}_{t \in [0,1]}$ be a ring motion of $H_2$ in~$\Rr^3 \setminus H_1$. 
We give $H_2$ the orientation induced from the $yz$-axis. 
We can give an orientation to $L_t$ which is induced from the orientation of~$H_2$.  
Let $c(L_t) \in \Rr^3$ be the center of~$L_t$, 
$r(L_t) \in \Rr_{>0}$ the radius, and $g^+(L_t)$ an element 
of the Grassman manifold $G^+(2,3)$ 
of oriented $2$-planes through the origin $O$ in~$\Rr^3$,
which is obtained from the oriented plane $H(L_t)$ containing $L_t$  
by sliding it along the vector~$\overrightarrow{c(L_t)O}$.  
Let $D(L_t)$ be the oriented disk in $\Rr^3$ bounded by $L_t$ 
in the plane $H(L_2)$ 
and let $d(L_t)$ be the intersection point~$D(L_t) \cap H_1$.  
Give $H_1$ an orientation induced from the $xy$-plane. 
Since each disk $D(L_t)$ intersects with $H_1$ on $d(L_t)$ in the positive direction, 
we can deform, up to homotopy as ring motions, 
the ring motion so that the normal vector of $D(L_t)$ at $d(L_t)$ is the positive unit tangent vector of~$H_1$.  
Then each $H(L_t)$ becomes a $2$-plane in $\Rr^3$ containing the $z$-axis.  
Next, we deform the ring motion so that the radius $r(L_t)$ is $1$ for all~$t \in [0,1]$.  
Finally, we deform the ring motion so that the center $c(L_t)$ is the intersection point~$d(L_t)$. 
Now we see that any ring motion is homotopic 
as ring motions to a normal ring motion. 
This implies that $R(\Rr^3 \setminus H_1, H_2)$ is generated by~$\ell$. 
\endproof

Now we discuss the ring group~$R(\Rr^3, H_1, H_2)$.  Let $H$ denote the Hopf link~$H_1 \sqcup H_2$. 
Let $\tau_H \in R(\Rr^3, H_1, H_2)$ be represented by~$\{ R_y(\pi t)(H) \}_{t \in [0,1]}$, 
\ie, the rotation of $180$ degrees about the $y$-axis 
and let $\ell \in R(\Rr^3, H_1, H_2)$ be represented by~$\{ R_z(2 \pi t)(H_2) \}_{t \in [0,1]}$, 
\ie, the positive standard rotation of $H_2$ along~$H_1$.  

\begin{lem} 
\label{L:Order}
In the ring group~$R(\Rr^3, H_1, H_2)$, we have the following. 
\begin{enumerate}[label={(\arabic*)}]
\item \label{I:first} $\tau_H^2 = \ell$ and~$\tau_H^4 = \ell^2 = 1$. 
\item \label{I:second} $\tau_H \ell \tau_H^{-1} = \ell^{-1}$.  
\item \label{I:third} The order of $\tau_H$ is $4$ and the order of $\ell$ is~$2$.  
\end{enumerate}
\end{lem} 


\proof 
\ref{I:first}
Let $f_{\tau_H}\colon [0,1] \to SO(3)$ 
and $f_{\ell}\colon [0,1] \to SO(3)$ be maps defined by 
\[
f_{\tau_H}(t) = R_y(\pi t) \in SO(3) \quad \mbox{and} \quad 
f_{\ell}(t) = R_z( 2 \pi t) \in SO(3).
\]
Then $[ f_{\tau_H} \ast f_{\tau_H} ] = [ f_{\ell}]  = -1 $ in~$\pi_1(SO(3)) = \{ 1, -1\}$. 
This implies that   $\tau_H^2 = \ell$ and $\tau_H^4 = \ell^2 = 1$ in~$R(\Rr^3, H_1, H_2)$. 

\ref{I:second} 
$[ f_{\tau_H}^{-1} \ast f_{\tau_H} \ast f_{\tau_H} ] = [ f_{\ell}^{-1}]  = -1 $ in~$\pi_1(SO(3)) = \{ 1, -1\}$. 

\ref{I:third} 
Consider the image of $\tau_H$ in the motion group $\M(\Rr^3, H)$ under the homomorphisms 
$R(\Rr^3, H_1, H_2) \to R(\Rr^3, H) \to \M(\Rr^3, H)$.  By using the double linking number
defined in~\cite{Carter-Kamada-Saito-Satoh:2002}, 
it can be seen that the order of $\tau_H$ in $\M(\Rr^3, H)$ is~$4$.
Thus the order of $\tau_4$ in $R(\Rr^3, H_1, H_2)$ is~$4$. By \ref{I:first}, the order of $\ell$ is $2$.  

\endproof

\begin{lem}
\label{L:ExactSequenceOrderedHopf}
Let $H_1$ and $H_2$ be the unit rings as above.  
Consider the composition of $e$ and~$p_1$:
\begin{equation}
\label{E:epOrderedHopf}
 \begin{CD}
R(\Rr^3 \setminus H_1, H_2)  @>e>>  R(\Rr^3, H_1, H_2) @>p_1>> R(\Rr^3, H_1).
\end{CD} 
\end{equation} 

\begin{enumerate}[label={(\arabic*)}]

\item \label{I:ExactFirst}
The sequence \eqref{E:epOrderedHopf} is exact. 

\item \label{I:ExactSecond}
The map $p_1$ is surjective. 
\end{enumerate}
\end{lem} 

\proof
\ref{I:ExactFirst} 
By Proposition~\ref{P:ExactSequenceRing}, we have that ${\mathrm Im}~{e} \subset {\mathrm Ker}~{p_1}$.  
We show that ${\mathrm Ker}~{p_1} \subset {\mathrm Im}~{e}$.  
Let $[ \{ L_{1(t)} \}_{t \in [0,1]} \sqcup \{ L_{2(t)} \}_{t \in [0,1]} ]$ belong to~${\mathrm Ker}~{p_1}$. 
Then $[ \{ L_{1(t)} \}_{t \in [0,1]}] =1$ in~$R(\Rr^3, H_1)$.  Thus the ring motion  $\{ L_{1(t)} \}_{t \in [0,1]}$ 
is homotopic to the stationary motion $\{ H_1 \}_{t \in [0,1]}$ of~$H_1$. To obtain such a homotopy, we use the strategy 
in the proof of Lemma~\ref{L:RingCircle}. Namely, first we change the ring motion $\{ L_{1(t)} \}_{t \in [0,1]}$ so that the center 
$c(L_{1(t)})$ of the ring $L_{1(t)}$ is the origin $O$ for every~$t \in [0,1]$, then change the radius~$r(L_{1(t)})$,
and change the element $g(L_{1(t)})$ of the Grassman manifold~$G(2,3)$. 
This procedure may change $\{ L_{2(t)} \}_{t \in [0,1]}$ by a homotopy 
keeping $L_{2(t)}$ to be a ring for every~$t$.  
Thus the ring motion $ \{ L_{1(t)} \}_{t \in [0,1]} \sqcup \{ L_{2(t)} \}_{t \in [0,1]} $ 
is equivalent to a motion which is the union of the statinary motion of $H_1$ and a ring motion of~$H_2$.  
Therefore,~${\mathrm Ker}~{p_1} \subset {\mathrm Im}~{e}$.  

\ref{I:ExactSecond}
By Lemma~\ref{L:RingCircle}, the ring group $R(\Rr^3, H_1)$ is generated by~$\tau_{H_1}$. 
The map $p_1$ sends 
$\tau_H \in R(\Rr^3, H_1, H_2)$ to~$\tau_{H_1} \in R(\Rr^3, H_1)$.  Thus $p_1$ is surjective. 
\endproof

\begin{thm} 
\label{T:PresentationOrderedHopf}
The ring group $R(\Rr^3, H_1, H_2)$ of the ordered Hopf link $H= H_1 \sqcup H_2$ admits the presentation 
\begin{equation} 
\label{E:OrderHopf}
\langle \tau_H \mid \tau_H^4 =1 \rangle. 
\end{equation}
\end{thm}

\proof 
By Lemma~\ref{L:ExactSequenceOrderedHopf}, we have a short exact sequence:  
\begin{equation} 
\label{E:OrderHopfA}
\begin{CD}
1 \longrightarrow e(R(\Rr^3 \setminus H_1, H_2)) 
  @> \iota >>  R(\Rr^3, H_1, H_2) @>p_1>> R(\Rr^3, H_1) 
  \longrightarrow 1. 
  \end{CD} 
\end{equation}
Since $R(\Rr^3 \setminus H_1, H_2)$ is generated by $\ell \in R(\Rr^3 \setminus H_1, H_2)$ (Lemma~\ref{L:MH2}), 
the image $e(R(\Rr^3 \setminus H_1, H_2))$ is generated by~$\ell \in R(\Rr^3, H_1, H_2)$.  
By Lemma~\ref{L:Order} the order of $\ell \in  R(\Rr^3, H_1, H_2)$ is~$2$.  Thus, we have 
\begin{equation} 
\label{E:OrderHopfB}
e(R(\Rr^3 \setminus H_1, H_2)) =  \langle \ell \mid \ell^2 =1 \rangle 
\end{equation} 
By Lemma~\ref{L:RingCircle}, we have 
\begin{equation} 
\label{E:OrderHopfC}
R(\Rr^3, H_1) = \langle \tau_{H_1} \mid \tau_{H_1}^2 = 1 \rangle. 
\end{equation} 
We choose $\tau_H$ as a set-theoretical lift of $\tau_{H_1}$ under~$p_1$.  By Lemma~\ref{L:Order}, we have 
\begin{equation}
\label{E:OrderHopfD}
\tau_H^2 = \ell \quad \mbox{and} \quad \tau_H \ell \tau_H^{-1} = \ell^{-1}.  
\end{equation}
Using the short exact sequence~\eqref{E:OrderHopfA}, presentations 
\eqref{E:OrderHopfB} and~\eqref{E:OrderHopfC}, and relations~\eqref{E:OrderHopfD}, 
and applying a standard method to give presentations for group extensions~\cite[Chapter~10]{Johnson:1997} we have that 
\begin{equation} 
\label{E:OrderHopfE}
R(\Rr^3, H_1, H_2) = \langle \ell, \tau_H \mid \ell^2 =1, \tau_H^2 = \ell, \tau_H \ell \tau_H^{-1} = \ell^{-1} \rangle, 
\end{equation}
which is reduced to the desired presentation~\eqref{E:OrderHopf}. 
\endproof 

Now we discuss the ring group~$R(\Rr^3, H)$.  Let $e_2$ be the unit vector $(0,1,0)$ in~$\Rr^3$.  

We consider an element $s \in R(\Rr^3, H)$ which interchanges $H_1$ and~$H_2$, 
represented by the ring motion  realized by a sequence of isometries of $\Rr^3$ as follows: 
first slide $H$ along~$(-1/2)e_2$, 
apply the rotation by $45$ degrees about the $y$-axis,
apply the rotation by $180$ degrees about the $x$-axis,
apply the rotation by $-45$ degrees about the $y$-axis, 
and slide along $(1/2)e_2$.
(This ring motion is equivalent to the following ring motion:   
first slide $H$ along $-e_2$, apply the rotation by $180$ degrees about the $x$-axis, 
and then apply the rotation by $-90$ degrees about the $y$-axis.) 

\begin{lem}
\label{L:PresentationHopf}
In the group $R(\Rr^3, H)$, 
the following relations are satisfied.
\begin{equation}
\label{E:HopfA}
s^2 = \tau_H^2 \quad \mbox{and} \quad 
s \tau_H s^{-1} = \tau_H^{-1} \quad 
\in R(\Rr^3, H). 
\end{equation}
\end{lem}  

\proof
We have $s^2 = \ell$ in~$R(\Rr^3, H)$
(it is easily understood when we draw link diagrams on the $yz$-plane). Thus,  by Lemma~\ref{L:Order}, 
we have~$s^2 = \tau_H^2$.  
By the sequence of isometries of $\Rr^3$ in the definition of~$s$, 
 the $y$-axis is sent to itself with reversed orientation. 
Since $\tau_H$ is a motion of $H$ realized by the rotation of $180$ degrees along the $y$-axis, we have~$s \tau_H s^{-1} = \tau_H^{-1}$.  
\endproof

\begin{thm}
\label{T:Hopfpresentation}
The ring group $R(\Rr^3, H)$ of the Hopf link admits the presentation 
\begin{equation}
\label{E:HopfB}
\langle \tau_H, s \mid \tau_H^4=1, ~ s^2 =\tau_H^2, ~ s \tau_H s^{-1} = \tau_H^{-1} \rangle. 
\end{equation}
\end{thm}

Remark that presentation \eqref{E:HopfB} can be rewritten as
\begin{equation}
\label{E:HopfBQuat}
\langle \tau_H, s \mid \tau_H^2  = s^2 = (\tau_H s)^2 \rangle, 
\end{equation}
which is a famous presentation of the \emph{quaternion group}.  

\proof 
The ring group $R(\Rr^3, H_1, H_2)$ is a subgroup of $R(\Rr^3, H)$ of index~$2$; 
consider the short exact sequence 
\begin{equation}
\label{E:Hopfexact}
1 \longrightarrow R(\Rr^3, H_1, H_2) \longrightarrow R(\Rr^3, H) \longrightarrow \Zz_2 \longrightarrow 1. 
\end{equation}
As a set-theoretical lift of the generator of $\Zz_2$ under  $R(\Rr^3, H) \to \Zz_2$, we choose $s \in R(\Rr^3, H)$.  
Using the short exact sequence \eqref{E:Hopfexact},  presentation \eqref{E:OrderHopf} of~$R(\Rr^3, H_1, H_2)$,
and relations \eqref{E:HopfA},  
we have  that $R(\Rr^3, H)$ admits the desired presentation~\eqref{E:HopfB}.  
\endproof

\section{The ring group of a Hopf link with a ring}
\label{S:HopfCircle}

In this section we discuss the ring group of an H-trivial link of type~$(1,1)$, \ie, the split union of a Hopf link and a ring. 

Let $H = H_1 \sqcup H_2$ be a Hopf link and $C$ a ring with~$H_1 = \{ (x, y, 0) \in \Rr^3 \mid x^2 + y^2 =1 \}$,  
$H_1 = \{ (0, y, z) \in \Rr^3 \mid (y-1)^2 + z^2 =1 \}$ and 
$C = \{ (x, y, 0) \in \Rr^3 \mid x^2 + (y-5)^2 =1\}$.   

\subsection{An exact sequence for  \texorpdfstring{$R(\Rr^3, H, C$)}{}}

\begin{lem}
\label{L:ExactSequenceHopfCircle}
Let $H$ and $C$ be as above.  
Consider the composition of $e$ and~$p_1$:
\begin{equation}
\label{E:epHopfCircle}
 \begin{CD}
R(\Rr^3 \setminus H, C)  @>e>>  R(\Rr^3, H, C) @>p_1>> R(\Rr^3, H).
\end{CD} 
\end{equation} 

\begin{enumerate}[label={(\arabic*)}]
\item \label{I:ExactHopfCircleFirst}
The sequence \eqref{E:epHopfCircle} is exact. 

\item \label{I:ExactHopfCircleSecond}
The map $p_1$ is surjective. 
\end{enumerate}
\end{lem}

\proof 
\ref{I:ExactHopfCircleFirst} 
By Proposition~\ref{P:ExactSequenceRing}, we have that ${\mathrm Im}~{e} \subset {\mathrm Ker}~{p_1}$.  
We show that ${\mathrm Ker}~{p_1} \subset {\mathrm Im}~{e}$.  
Let $[ \{ L_{1(t)} \}_{t \in [0,1]} \sqcup \{ L_{2(t)} \}_{t \in [0,1]} ]$ belong to ${\mathrm Ker}~{p_1}$. 
Then $[ \{ L_{1(t)} \}_{t \in [0,1]}] =1$ in $R(\Rr^3, H)$.  Thus the ring motion  $\{ L_{1(t)} \}_{t \in [0,1]}$ 
is homotopic to the stationary motion $\{ H \}_{t \in [0,1]}$ of $H$. 
We show that $\{ L_{1(t)} \}_{t \in [0,1]} \sqcup \{ L_{2(t)} \}_{t \in [0,1]}$ are homotopic as ring motions to 
the union of the stationary motion of $H$ and a ring motion of $C$.

\begin{step}
First deform the ring motion  $\{ L_{1(t)} \}_{t \in [0,1]}$ of $H$ and the motion $\{ L_{2(t)} \}_{t \in [0,1]}$ of $C$ in such a way that the restriction to $H_1$ becomes a stationary motion of $H_1$ keeping the condition that 
the new $\{ L_{1(t)} \}_{t \in [0,1]}$ and $\{ L_{2(t)} \}_{t \in [0,1]}$ are disjoint ring motions.  This is done by the strategy 
used
in the proof of Lemma~\ref{L:RingCircle} to deform the motion of $H_1$ to the stationary motion. (Recall the proof of Lemma~\ref{L:ExactSequenceOrderedHopf}.)  

Now we may assume that the restriction of $\{ L_{1(t)} \}_{t \in [0,1]}$ to $H_1$ is the stationary motion.  
The restriction of $\{ L_{1(t)} \}_{t \in [0,1]}$ to $H_2$ is a ring motion of $H_2$ in~$\Rr^3 \setminus H_1$.  
\end{step}

\begin{step}
Secondly, deform the ring motion  $\{ L_{1(t)} \}_{t \in [0,1]}$ of $H$ and the motion $\{ L_{2(t)} \}_{t \in [0,1]}$ of $C$ so that the restriction to $H$ becomes the stationary motion of $H$ keeping the condition that the new $\{ L_{1(t)} \}_{t \in [0,1]}$ and $\{ L_{2(t)} \}_{t \in [0,1]}$ are disjoint ring motions.  This is done as follows: Since the restriction of $\{ L_{1(t)} \}_{t \in [0,1]}$ to $H_2$ is a ring motion of $H_2$ in~$\Rr^3 \setminus H_1$, it is homotopic to the power of the positive or negative standard rotation of $H_2$ along $H_1$ by the argument in the proof of Lemma~\ref{L:MH2}.   During the homotopy for~$\{ L_{1(t)} \}_{t \in [0,1]}$, we may deform 
$\{ L_{2(t)} \}_{t \in [0,1]}$ keeping the condition that it is a ring motion.  

Now, $\{ L_{1(t)} \}_{t \in [0,1]} \sqcup \{ L_{2(t)} \}_{t \in [0,1]}$ satisfies that $\{ L_{1(t)} \}_{t \in [0,1]}$ is the stationary motion of $H$ and $\{ L_{2(t)} \}_{t \in [0,1]}$ is a ring motion of $H_2$ in $\Rr^3 \setminus H$. Thus it represents an element of the image of $e: R(\Rr^3 \setminus H, C)  \to   R(\Rr^3, H, C)$.  

\ref{I:ExactHopfCircleSecond} 
By Lemma~\ref{T:Hopfpresentation}, the ring group $R(\Rr^3, H)$ is generated by $\tau_H$ and~$s$. 

Let $\tilde\tau_H$ (or~$\tilde s$)  
be elements of $R(\Rr^3, H, C)$ which is the union of $\tau_H$ (or $s$) and the stationary motion on~$C$.   
Then $p_1 (\tilde\tau_H) = \tau_H$ and~$p_1(\tilde s) =s$. 
Thus $p_1$ is surjective. 

\end{step}
\endproof 

Later, in Lemma~\ref{L:ExactSequenceHopfCircleB}, we will see 
that sequence \eqref{E:epHopfCircle} induces a short exact sequence  
that will allow us to use once more the standard method to write presentations of group extensions.

\subsection{On the ring group  \texorpdfstring{$R(\Rr^3 \setminus H, C)$}{}}

\begin{figure}[ht]
\centering
\includegraphics[scale=.8]{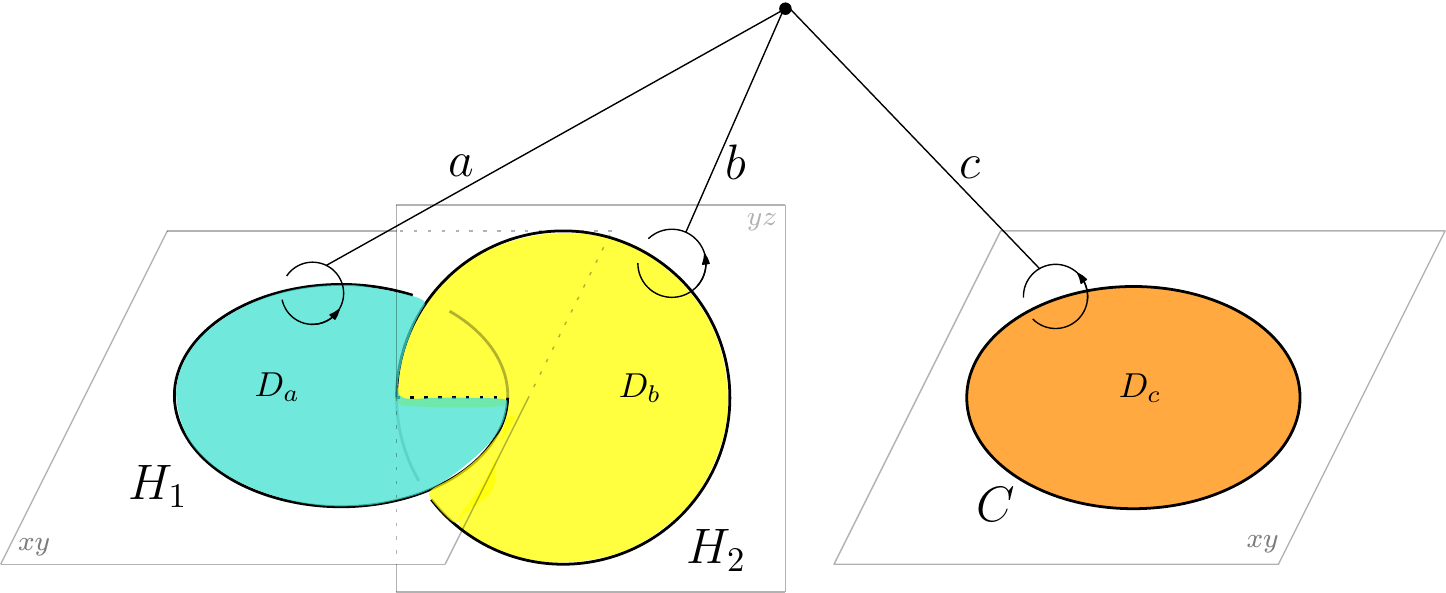}
\caption{Generators of~$\pi_1(\Rr^3 \setminus (H \sqcup C))$. }
\label{F:HunionC}
\end{figure}

Let $H= H_1 \sqcup H_2$ be the Hopf link and $C$ the ring disjoint from~$H$ as before.
Let us choose a base point for the fundamental group~$\pi_1(\Rr^3 \setminus (H \sqcup C))$ in such a way that the $z$-coordinate is sufficiently large. 
Let~$a$, $b$ and $c$ be elements of $\pi_1(\Rr^3 \setminus (H \sqcup C))$ represented by meridian loops of~$H_1$,~$H_2$, and~$C$, respectively. We assume that these meridian loops are oriented such that the linking number is $+1$ when we give~$H_1$,~$H_2$, and $C$
orientations induced from the $xy$-plane and the $yz$-plane, see Figure~\ref{F:HunionC}.  
The fundamental group is the free product of the free abelian group of rank $2$ generated by $a$ and $b$ and the infinite cyclic group generated by~$c$: 
\[
\pi_1(\Rr^3 \setminus (H \sqcup C)) = 
\langle a, b, c \mid [a, b]=1 \rangle \cong (\Zz \oplus \Zz) \ast \Zz. 
\]

\begin{figure}[hbtp]
\centering
\includegraphics[scale=.5]{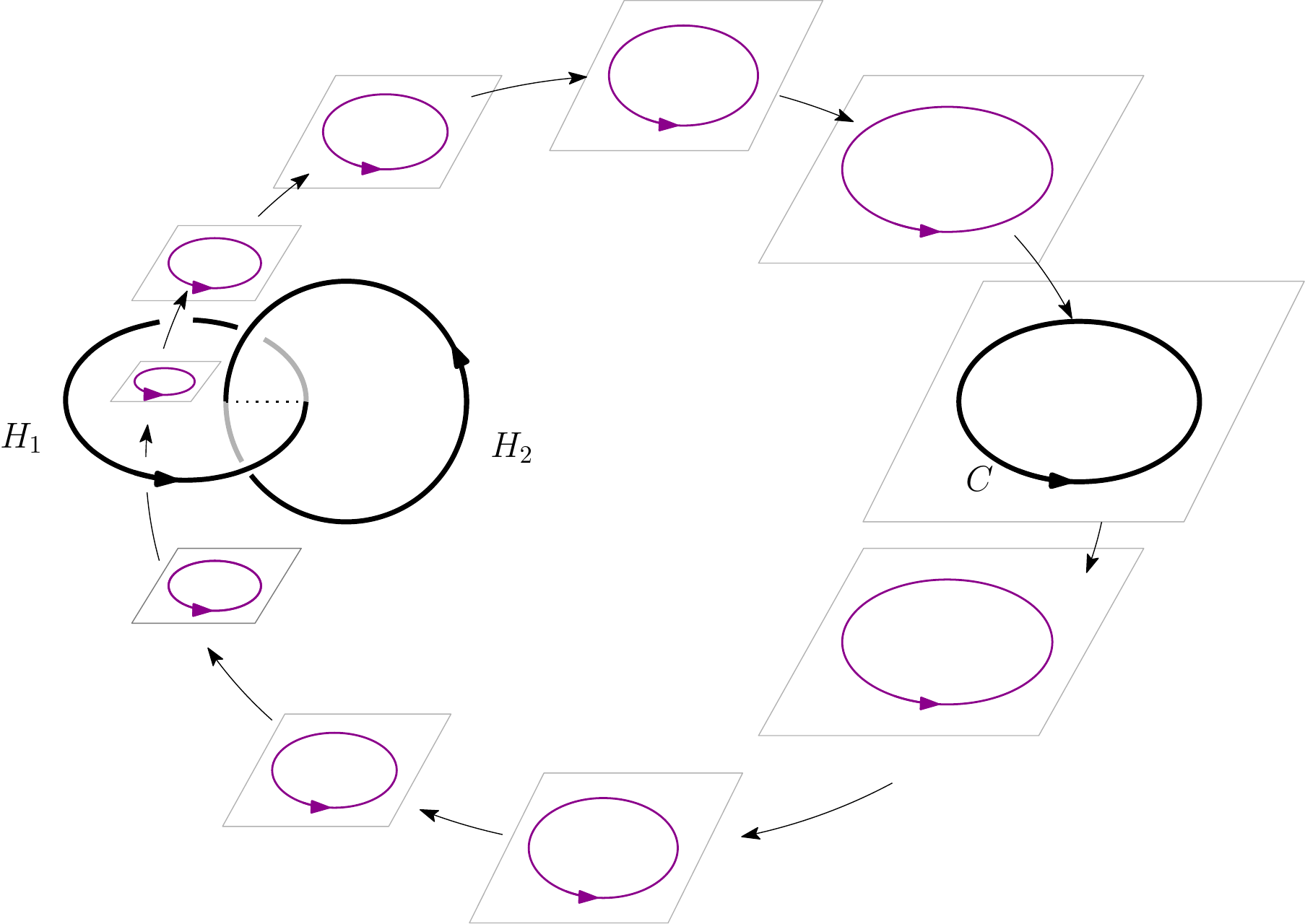}
\caption{The ring motion~$g_a$. }
\label{F:g_a}
\end{figure}

\begin{figure}[hbtp]
\centering
\includegraphics[scale=.5]{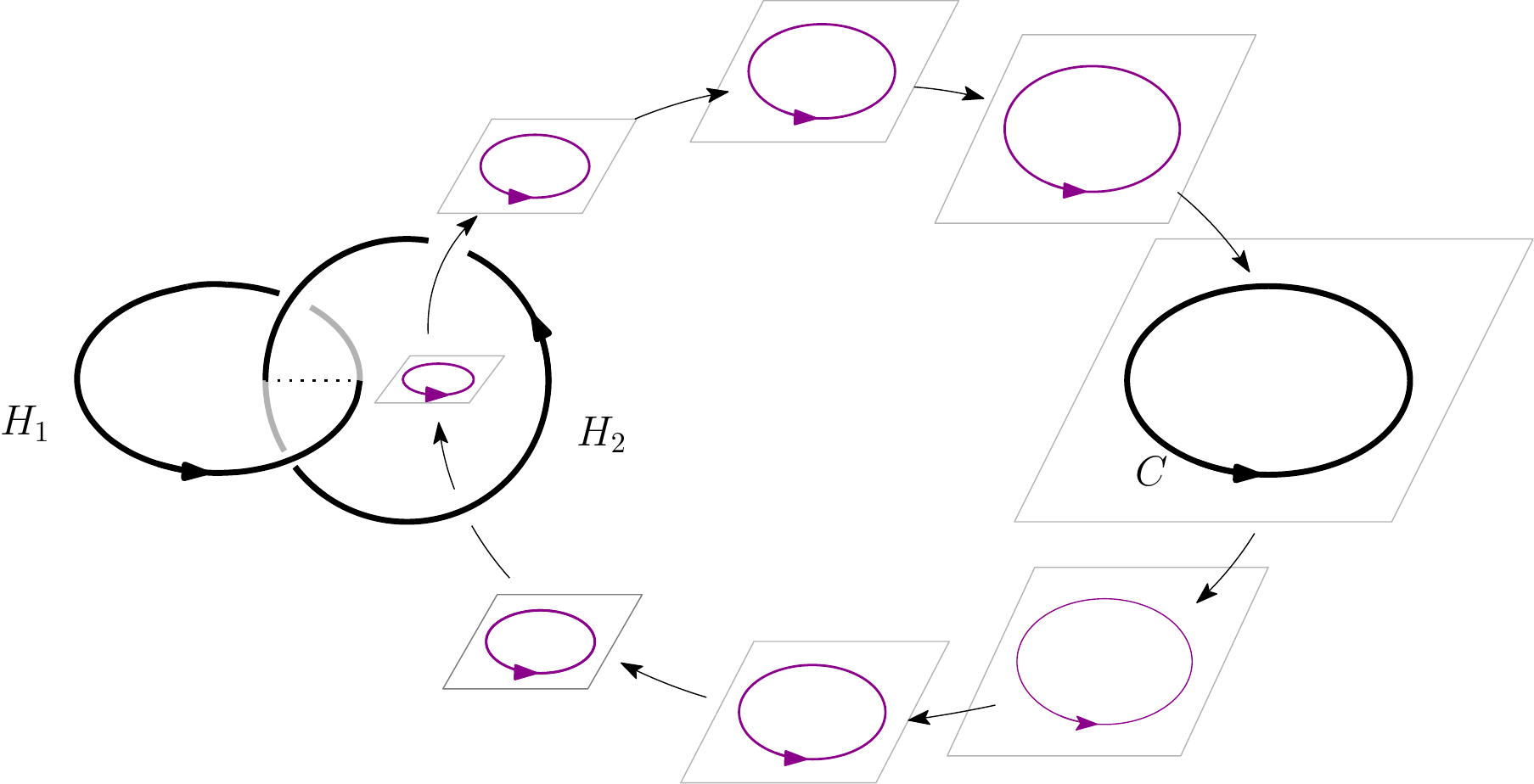}
\caption{The ring motion~$g_b$. }
\label{F:g_b}
\end{figure}

\begin{figure}[hbtp]
\centering
\includegraphics[scale=.4]{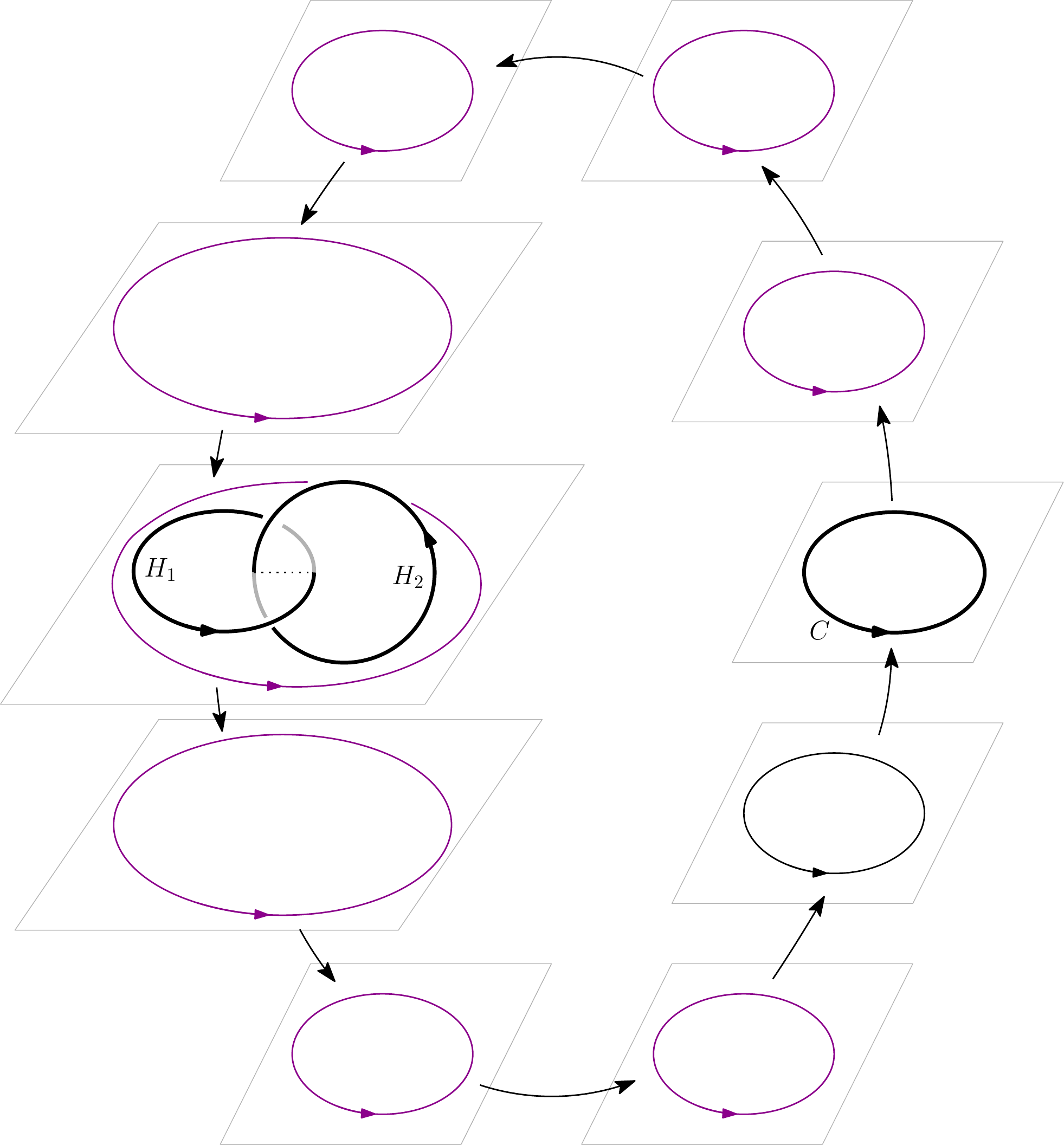}
\caption{The ring motion~$\varepsilon_C$. }
\label{F:epsilon_C}
\end{figure}

Let us introduce some ring motions:
\begin{itemize}
\item $g_a$: $C$ pulls through~$H_1$, see Figure~\ref{F:g_a};
\item $g_b$: $C$ pulls through~$H_2$, see Figure~\ref{F:g_b};
\item $\tau_C$: $C$ rotates by~$180$ degrees around the $y$-axis, as in Section~\ref{S:Circle};
\item $\varepsilon_C$: $C$ traslates above~$H$, slides downwards encircling~$H$,
and then traslates back to its original position, see Figure~\ref{F:epsilon_C}.
\end{itemize}

\begin{lem} 
\label{L:RingGroupFirstHCgenerator} 
The ring group $R(\Rr^3 \setminus H, C)$ is generated by 
$g_a, g_b, \varepsilon_{C}, \tau_C$. 
The following relations are satisfied. 
\begin{equation}
 [g_a, g_b] = 1, \tau_C^2=1, \, 
 [g_a, \tau_C] = 
 [g_b, \tau_C] =  1, 
 \tau_C \varepsilon_{C} \tau_C = \varepsilon_{C}^{-1}.  
\end{equation}
\end{lem}

\proof  
First of all, remark that the motion group $R(\Rr^3 \setminus H, \ast)$, 
where $\ast$ is a point,
is the fundamental group $\pi_1(\Rr^3 \setminus H) = \langle a, b \mid [a, b]=1 \rangle \cong (\Zz \oplus \Zz) $, 
and recall that $R(\Rr^3, C) = \langle \tau_C \mid \tau_C^2=1 \rangle \cong \Zz_2$ (Section~\ref{S:Circle}).

Consider~$D_a$,~$D_b$ and $D_c$ to be disks bounded by~$H_1$, $H_2$ and~$C$, 
flatly embedded in the planes where $H_1$, $H_2$ and~$C$ lie, 
as in Figure~\ref{F:HunionC}.
Let  $\{ C_t \}_{t \in [0,1]}$ be a ring motion of $C$ in~$\Rr^3 \setminus H$, 
and $D_{C_t}$ be the flat disk bounded by~$C_t$, for~$t \in [0, 1]$.
Let us distinguish two cases.

\begin{enumerate}[label={\arabic*.}]
\item Suppose that for all~$t \in [0, 1]$, 
$D_{C_t} \cap (D_a \cup D_b) = \emptyset$. 
After a deformation of~$\{ C_t \}_{t \in [0,1]}$ by a homotopy,
we may assume that there exists a convex $3$-ball~$B_C$, 
disjoint from $(D_a \cup D_b)$, 
and such that $C_t$ lies in $B_C$ for all~$t \in [0, 1]$.
Then 
$\{ C_t \}_{t \in [0,1]}$ 
represents 
an element of $R(B_C, C) \cong R(\Rr^3, C) = \langle \tau_C \mid \tau_C^2=1 \rangle$.

\item \label{I:intersection} Suppose that for some value of~$t$, 
$D_{C_t} \cap (D_a \cup D_b) \neq \emptyset$. Then let us consider these two subcases.

\begin{enumerate}[label={\ref{I:intersection}\arabic*.}]
\item The disks $D_{C_t}$ intersects the interior of $D_a$ and/or $D_b$ for $t$ 
in a finite number of intervals~$[\tilde{t}-\varepsilon, \tilde{t}+\varepsilon]$, 
and $H_1 \cap \mathrm{int}(D_{C_t}) = H_2 \cap \mathrm{int}(D_{C_t}) = \emptyset$ 
for all~$t \in [0, 1]$. Then  
$\{ C_t \}_{t \in [0,1]}$, modulo $\tau_C$, represents 
an element of~$R(\Rr^3 \setminus H, \ast)=\langle a, b \mid [a, b]=1 \rangle$.

\item The interiors $\mathrm{int}(D_{C_t})$ intersects $H_1$ and/or $H_2$ for $t$ 
in a finite number of intervals~$[\tilde{t}-\varepsilon, \tilde{t}+\varepsilon]$, 
and $C_t \cap (\mathrm{int}(D_a) \cup \mathrm{int}(D_b) )= \emptyset$ 
for all~$t \in [0, 1]$. Then 
$\{ C_t \}_{t \in [0,1]}$, modulo~$\tau_C$, represents
an element of the subgroup of $R(\Rr^3 \setminus H, C)$ 
generated by the motion~$\varepsilon_C$ (Figure~\ref{F:epsilon_C}). 

\end{enumerate}
\end{enumerate}

Every generic ring motion of $C$ in $\Rr^3 \setminus H$ can be decomposed
in a combination of motions that fall in the considered cases, thus 
$\tau_C, g_a, g_b$ and $\varepsilon_C$ are a generating set for~$R(\Rr^3 \setminus H, C)$.

\begin{figure}[hbtp]
\centering
\includegraphics[scale=.5]{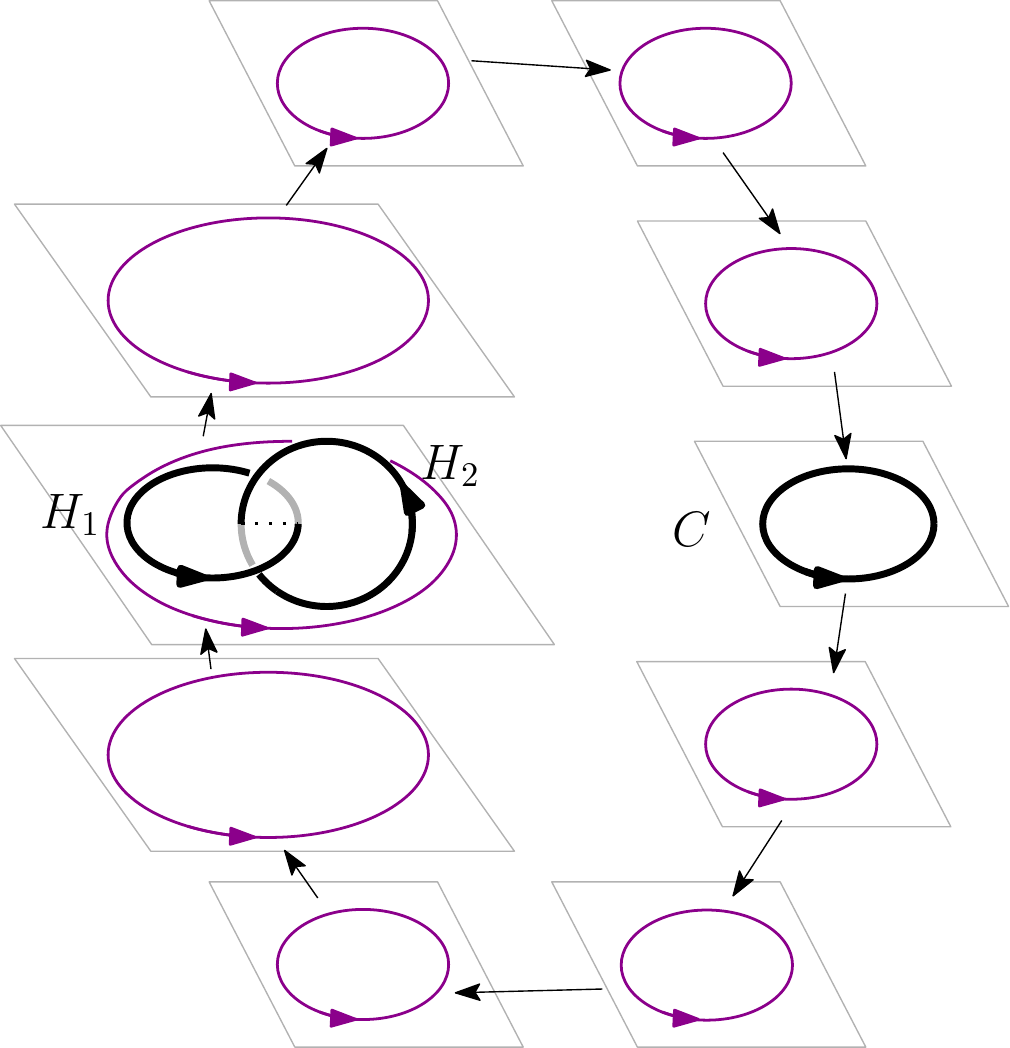}
\caption{The ring motion $\varepsilon_{C}^{-1}$. }
\label{F:Tau_Epsilon_Tau_1}
\end{figure}

\begin{figure}[hbtp]
\centering
\includegraphics[scale=.5]{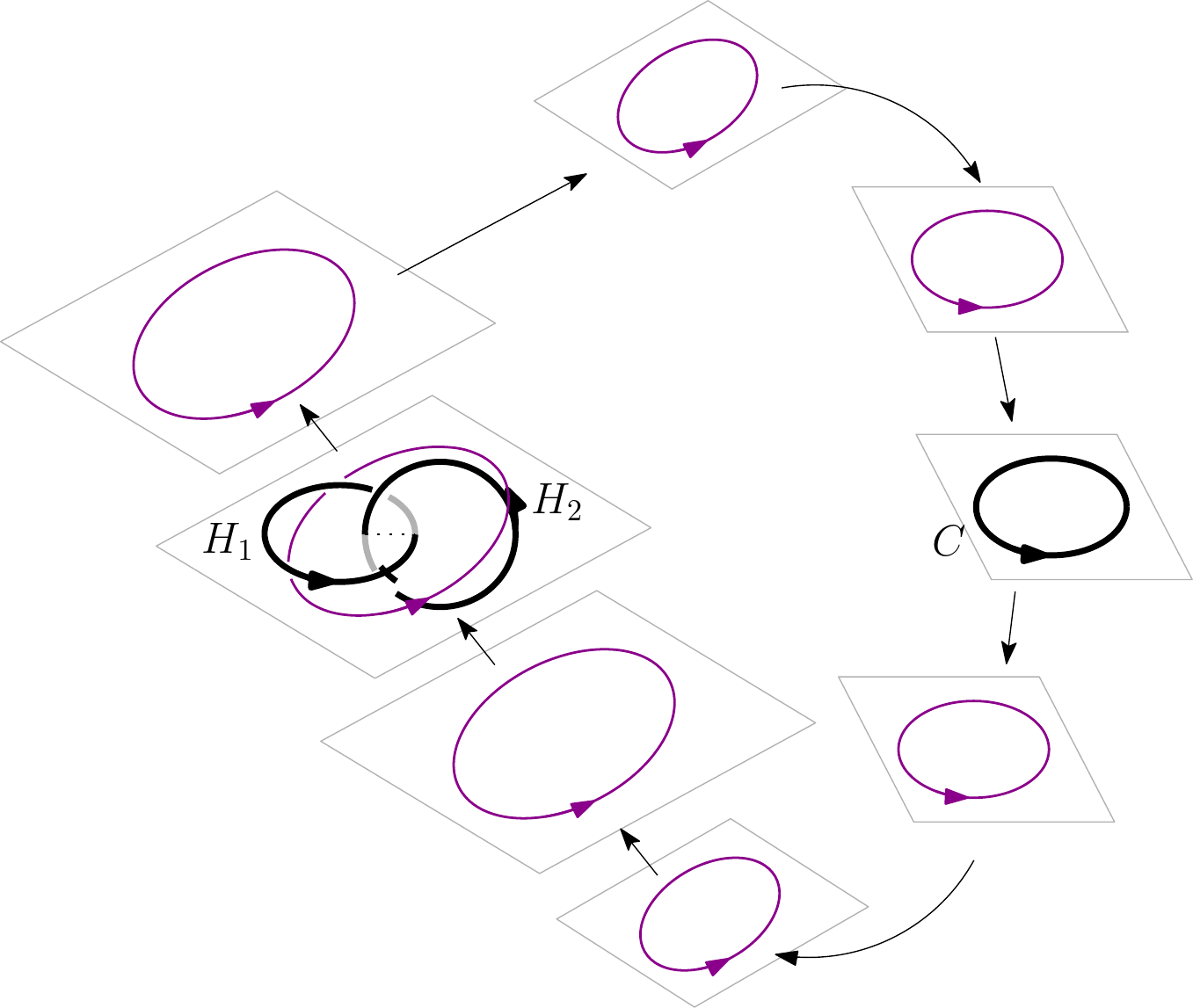}
\caption{A deformation of $\varepsilon_{C}^{-1}$, where the plane where $C$ is lying slightly tilts before encircling~$H$. }
\label{F:Tau_Epsilon_Tau_2}
\end{figure}

\begin{figure}[hbtp]
\centering
\includegraphics[scale=.5]{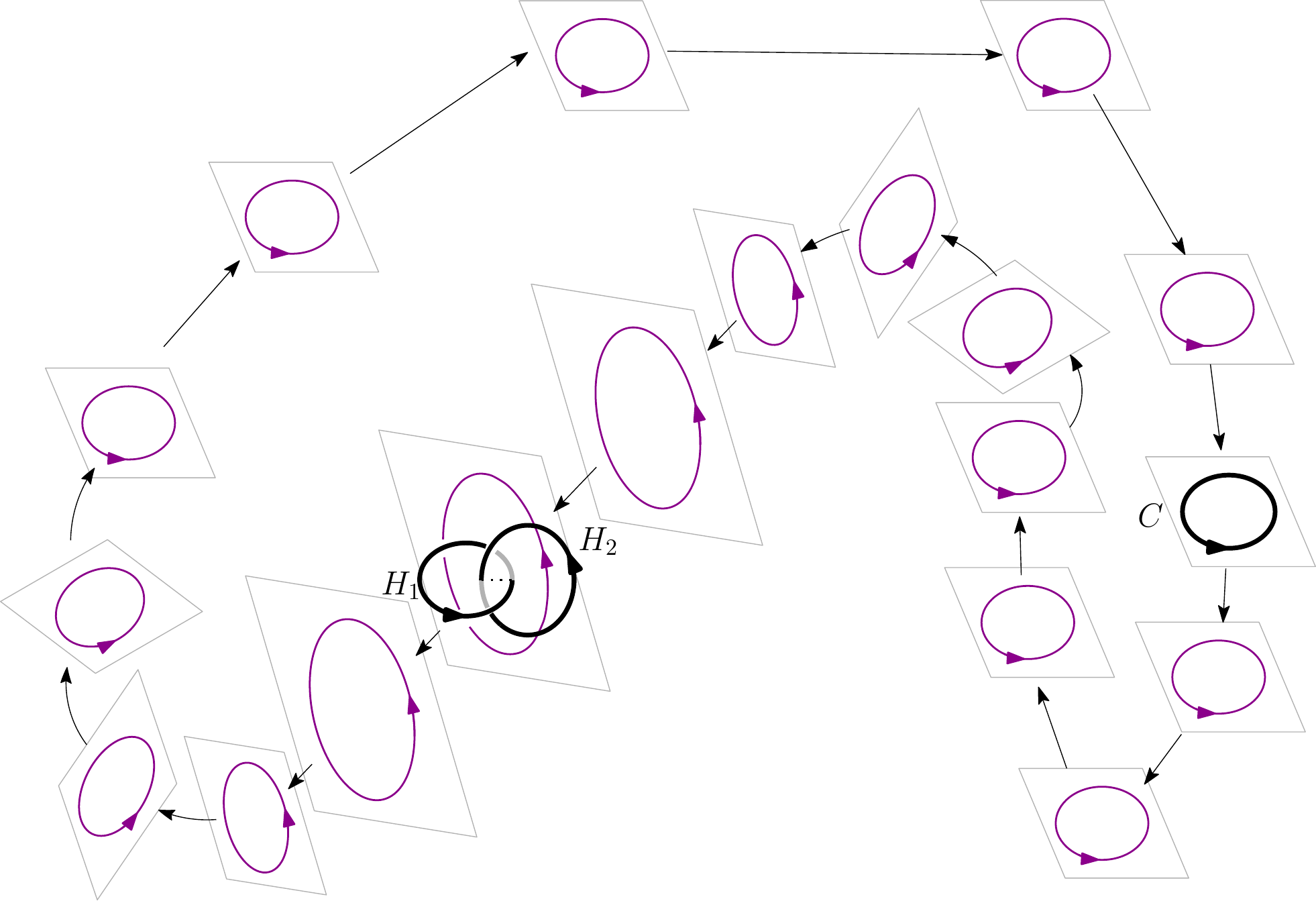}
\caption{The plane where $C$ is lying tilts a bit more before encircling~$H$. }
\label{F:Tau_Epsilon_Tau_3}
\end{figure}

\begin{figure}[hbtp]
\centering
\includegraphics[scale=.5]{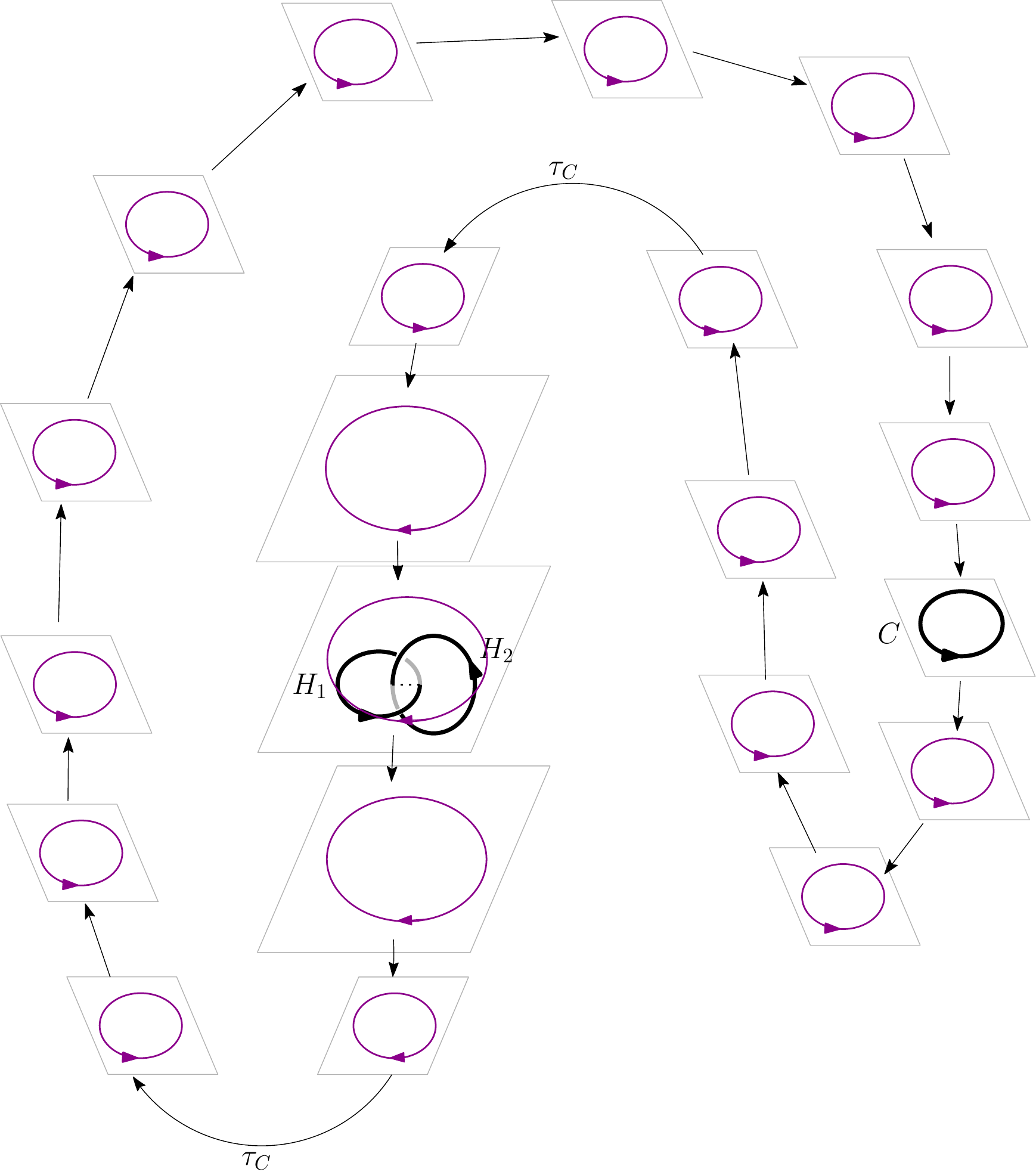}
\caption{ The plane where $C$ is lying tilts by $180$ degrees before encircling~$H$, and the ring motion $\varepsilon_{C}^{-1}$
has been continuously deformed into he motion~$\tau_C \varepsilon_{C} \tau_C $.}
\label{F:Tau_Epsilon_Tau_4}
\end{figure}

The relations in the statement descend from relations of  $R(\Rr^3 \setminus H, \ast)$  and~$R(\Rr^3, C)$, 
with the exception of~$\tau_C \varepsilon_{C} \tau_C = \varepsilon_{C}^{-1}$. 
This last relation can be seen from the sequence of Figures~\ref{F:Tau_Epsilon_Tau_1}, 
\ref{F:Tau_Epsilon_Tau_2}, \ref{F:Tau_Epsilon_Tau_3}, and~\ref{F:Tau_Epsilon_Tau_4}.
\endproof

Let $R^+(\Rr^3 \setminus H, C)$ be the index $2$ subgroup 
of $R(\Rr^3 \setminus H, C)$ 
consisting of equivalence classes of ring motions of $C$ 
that preserve an orientation of~$C$. 
This is the subgroup generated by~$g_a$,~$g_b$,~$\varepsilon_{C}$.

\begin{lem} 
\label{L:RingGroupFirstHCpresentation} 
The ring group $R^+(\Rr^3 \setminus H, C)$ admits the  presentation 
\begin{equation}
\label{E:Presentation-Hplus}
\langle \, g_a, g_b, \varepsilon_{C} \mid 
 [g_a, g_b] = 1
\, \rangle,   
\end{equation}
and the Dahm homomorphism 
$D \colon R^+(\Rr^3 \setminus H, C) \to {\mathrm Aut}(\pi_1(\Rr^3 \setminus (H \sqcup C)))$ is injective. 
\end{lem}

\proof 
The images of the elements $g_a$, $g_b$ and $\varepsilon_{C}$ under 
the Dahm homomorphism $D: R^+(\Rr^3 \setminus H, C) \to 
{\mathrm Aut}(\pi_1(\Rr^3 \setminus (H \sqcup C)))= {\mathrm Aut}\big(\langle a, b, c \mid [a, b]=1 \rangle\big)$ are the following automorphisms:
\begin{equation}
D(g_a) \colon
\left\{
\begin{array}{lll}
a & \mapsto & a \\ 
b & \mapsto & b \\
c & \mapsto & a c a^{-1} 
\end{array}
\right. 
\\ 
D(g_b) \colon
\left\{
\begin{array}{lll}
a & \mapsto & a \\ 
b & \mapsto & b \\
c & \mapsto & b c b^{-1} 
\end{array}
\right. 
\\ 
D(\varepsilon_{C}) \colon
\left\{
\begin{array}{lll}
a & \mapsto & c a c^{-1}\\ 
b & \mapsto & c b c^{-1}\\
c & \mapsto & c.  
\end{array}
\right. 
\end{equation}  

Let $G_1$ be the free abelian group generated by $g_1$ and~$g_2$,
let $G_2$ be the infinite cyclic group generated by~$\varepsilon_{C}$, 
and let $G$ be the free product of $G_1$ and~$G_2$, \ie, $G= \langle g_a, g_b, \varepsilon_{C} \mid [g_a, g_b]=1 \rangle$.  
We show that the natural epimorphism $\mu: G \to R^+(\Rr^3 \setminus H, C)$ is injective by showing that 
the homomorphism $D' = D \circ \mu : G \to {\mathrm Aut}(\pi_1(\Rr^3 \setminus (H \sqcup C)))$ is injective. 

Let $W\colon G \to \langle a, b, c \mid [a, b]=1 \rangle$ be the isomorphism with $g_a \mapsto a, g_b \mapsto b, \varepsilon_{C} \mapsto c$.  
Note that for any  $g \in G$, 
$D'(g)$ is the inner automorphism of $\langle a, b, c \mid [a, b]=1 \rangle$ by $W(g)$, \ie,
$D'(g)(x) = W(g) x W(g)^{-1}$. This implies that  
$D'(g)=1$ if and only if $W(g)=1$.  Thus, $D'$ is an isomorphism and we have the presentation~\eqref{E:Presentation-Hplus}. 
\endproof

\begin{rmk}
Remark that $\pi_1(\Rr^3 \setminus (H \sqcup C))$ is a  right-angled Artin group, and that 
$\{D(g_a), D(g_b), D(\varepsilon_C)\}$ is the set of (partial) conjugations in~${\mathrm Aut}(\pi_1(\Rr^3 \setminus (H \sqcup C)))$.
Then $\{D(g_a), D(g_b), D(\varepsilon_C)\}$ is a generating set for a particular case of \emph{group of vertex-conjugating automorphisms of a right-angled Artin group}, 
for which Toinet gives a complete presentation in~\cite{Toinet:2012}. In this paper he generalises a method used by McCool~\cite{McCool:1986} 
to study \emph{groups of basis-conjugating automorphisms of free groups}.
We recall that these last ones are isomorphic to \emph{pure untwisted} ring groups, and to \emph{pure} loop braid groups~\cite{BrendleHatcher:2013, Damiani:Journey}
\end{rmk}

\begin{lem} 
\label{L:RingGroupHCpresentation} 
The ring group $R(\Rr^3 \setminus H, C)$ admits the presentation
\begin{equation}
\label{E:Presentation-H}
\langle \, g_a, g_b, \varepsilon_{C}, \tau_{C} \mid 
 [g_a, g_b] = 1, \tau_{C}^2=1, \, 
 [g_a, \tau_{C}] = 
 [g_b, \tau_{C}] =  1, 
 \tau_{C} \varepsilon_{C} \tau_{C} = \varepsilon_{C}^{-1} 
\, \rangle,   
\end{equation}
and the Dahm homomorhism 
$D\colon R(\Rr^3 \setminus H, C) \to {\mathrm Aut}(\pi_1(\Rr^3 \setminus (H \sqcup C)))$ is injective. 
\end{lem}

\proof 
Presentation \eqref{E:Presentation-H} 
is obtained from presentation \eqref{E:Presentation-Hplus} and Lemma~\ref{L:RingGroupFirstHCgenerator} 
by using the short exact sequence 
\begin{equation} 
\label{E:seHC}
\begin{CD}
1 \longrightarrow R^+(\Rr^3 \setminus H, C) 
  @> \iota >>  R(\Rr^3 \setminus H, C) 
  @>>> \Zz_2
  \longrightarrow 1. 
  \end{CD} 
\end{equation}

Let $g \in R(\Rr^3 \setminus H, C)$ be an element of the kernel of~$D$. 
In Lemma~\ref{L:RingGroupFirstHCpresentation} we have seen that $D$ is injective
on the subgroup~$R^+(\Rr^3 \setminus H, C)$. 
Suppose $g \in R(\Rr^3 \setminus H, C) \setminus R^+(\Rr^3 \setminus H, C)$.  
Then $g = g_0 \tau_{C}$ for some~$g_0 \in R^+(\Rr^3 \setminus H, C)$.  Since 
\begin{equation}
D(\tau_{C}) \colon
\left\{
\begin{array}{lll}
a & \mapsto & a \\ 
b & \mapsto & b \\
c & \mapsto & c^{-1},  
\end{array}
\right. 
\end{equation}
$D(\tau_{C})$ is never an inner automorphism of~$\pi_1(\Rr^3 \setminus (H \sqcup C)))$. 
This contradicts to that $D(g_0)$ is an inner automorphism.  Thus, $D$ is injection on~$R(\Rr^3 \setminus H, C)$.  
\endproof

\begin{lem}
\label{L:ExactSequenceHopfCircleB}
The sequence involving $e$ and $p_1$ in Lemma~\ref{L:ExactSequenceHopfCircle} induces the short exact sequence 
\begin{equation}
\label{E:epHopfCircleB}
 \begin{CD}
 1 \longrightarrow 
R(\Rr^3 \setminus H, C)  @>e>>  R(\Rr^3, H, C) @>p_1>> R(\Rr^3, H) 
\longrightarrow 1.
\end{CD} 
\end{equation} 
\end{lem} 

\proof 
By Lemma~\ref{L:ExactSequenceHopfCircle}, it is sufficient to show that $e$ is injective.  
This  follows from the injectivity of the Dahm homomorphism 
$D\colon R(\Rr^3 \setminus H, C) \to \pi_1(\Rr^3 \setminus (H \sqcup C)))$.  
\endproof

\subsection{The ring group  \texorpdfstring{$R(\Rr^3, H \sqcup C)$}{}}

\begin{thm}
\label{T:HopfCirclepresentation}
The ring group $R(\Rr^3, H \sqcup C)$ $(= R(\Rr^3, H, C))$ admits the following presentation: 
Generators: 
\begin{equation}
\label{E:HCgenerator}
g_a, g_b, \varepsilon_C, \tau_C, \tau_H, s. 
\end{equation}
Relations: 
\begin{equation}
\label{E:HCrelationA}
 [g_a, g_b] = 1, ~ \tau_{C}^2=1, ~  
 [g_a, \tau_{C}] = 
 [g_b, \tau_{C}] =  1, ~ 
 \tau_{C} \varepsilon_{C} \tau_{C} = \varepsilon_{C}, 
 \end{equation}
 \begin{equation} 
 \label{E:HCrelationB}
 \tau_H^4 =1, ~ s^2 =\tau_H^2, ~ s \tau_H s^{-1} = \tau_H^{-1}, 
  \end{equation}
 \begin{equation} 
 \label{E:HCrelationC}
 \tau_H g_a \tau_H^{-1} = g_a^{-1}, ~ 
 \tau_H g_b \tau_H^{-1} = g_b^{-1}, ~ 
 \tau_H \varepsilon_C \tau_H^{-1} = \varepsilon_C, ~ 
 \tau_H \tau_C \tau_H^{-1} = \tau_C, 
  \end{equation}
 \begin{equation}  
 \label{E:HCrelationD}
 s g_a s^{-1} = g_a, ~ 
 s g_b s^{-1} = g_b, ~ 
 s \varepsilon_C s^{-1} = \varepsilon_C, ~ 
 s \tau_C s^{-1} = \tau_C.  
\end{equation}
\end{thm}

\proof 
Consider the short exact sequence \eqref{E:seHC}.  Let $\tilde\tau_H$ (or~$\tilde s$)  
be elements of $R(\Rr^3, H \sqcup C)$ which is the union of $\tau_H$ (or $s$) and the stationary motion on~$C$.   
Then $p_1 (\tilde\tau_H) = \tau_H$ and~$p_1(\tilde s) =s$.  We have a section 
$R(\Rr^3, H) \to R(\Rr^3, H \sqcup C)$ sending $\tau_H$ to $\tilde\tau_H$ and $s$ to~$\tilde s$.  
Thus, the short exact sequence \eqref{E:seHC} is split. We may denote the elements $\tilde\tau_H$ and $\tilde s$ by $\tau_H$ and $s$ for simplicity.   
Using the presentation \eqref{E:Presentation-H} of~$R(\Rr^3 \setminus H, C)$, and the presentation \eqref{E:HopfB} of~$R(\Rr^3, H)$, 
we have the generators~\eqref{E:HCgenerator} and relations \eqref{E:HCrelationA} and~\eqref{E:HCrelationB}. 
The actions of $\tau_H$ and $s$ yield relations \eqref{E:HCrelationC} and~\eqref{E:HCrelationD}.  
\endproof


\section*{Acknowledgements}
During the writing of this paper both authors were supported by JSPS KAKENHI Grant Number JP16F1679.
The first author was also supported by a JSPS Postdoctoral Fellowship For Foreign Researchers, and the second author was also supported by JSPS KAKENHI Grant Number JP26287013.
We thank Riccardo Piergallini and John Guaschi for helpful discussions, Eric Rowell for giving us access to a precious reference, and  Arnaud Mortier 
for the interest he expressed in this work. 

\bibliography{H_trivial.bib}{}
\bibliographystyle{alpha}

\end{document}